\documentclass{amsart}
\usepackage[english]{babel}
\usepackage[latin1]{inputenc}
\usepackage{mathrsfs}
\usepackage{amsmath}
\usepackage{amsthm}
\usepackage{amssymb}
\usepackage{indentfirst}

\usepackage[dvips]{graphicx}

\usepackage{ifthen}
\newcommand{\cit}[1]{{\rm \textbf{#1}}}

\newcommand{\Ref}[2]{\cit{%
\ifthenelse{\equal{#1}{thm}}{Theorem}{}%
\ifthenelse{\equal{#1}{prop}}{Proposition}{}%
\ifthenelse{\equal{#1}{lem}}{Lemma}{}%
\ifthenelse{\equal{#1}{cor}}{Corollary}{}%
\ifthenelse{\equal{#1}{defn}}{Definition}{}%
\ifthenelse{\equal{#1}{oss}}{Remark}{}%
\ifthenelse{\equal{#1}{sec}}{Section}{}%
\ifthenelse{\equal{#1}{ex}}{Example}{}%
\ifthenelse{\equal{#1}{conj}}{Conjecture}{}%
\ifthenelse{\equal{#1}{ssec}}{Subsection}{}%
\ifthenelse{\equal{#1}{tab}}{Table}{}%
\ifthenelse{\equal{#1}{cla}}{Claim}{}%
\  \ref{#1:#2}%
}}

\theoremstyle{plain} 
\newtheorem{prop}{Proposition}[section]
\newtheorem{thm}[prop]{Theorem}

\theoremstyle{remark}

\newtheorem{oss}[prop]{Remark}
\newtheorem{ex}[prop]{Example}

\theoremstyle{definition}
\newtheorem{defn}[prop]{Definition}
\newcommand{\hk}{hyperk\"{a}hler }

\newcommand{\kahl}{K\"{a}hler }

\newcommand{\kntiposp}{$K3^{[n]}$ type }
\newcommand{\ie}{i.~e.~}

\newcommand{\issp}{irreducible symplectic }

\begin{document}
\title{On the monodromy of irreducible symplectic manifolds}
\author{Giovanni Mongardi}
\address{Dipartimento di Matematica, Universit\'{a} degli studi di Milano,  Via Cesare Saldini 50, Milano, 20133 Italia  }
\begin{abstract}
Exploiting recent results on the ample cone of irreducible symplectic manifolds, we provide a different point of view for the computation of their monodromy groups. In particular, we give the final step in the computation of the monodromy group for generalised Kummer manifolds and we prove that the monodromy of O'Grady's ten dimensional manifold is smaller than what was expected.
\end{abstract}
\thanks{Supported by FIRB 2012 ``Spazi di Moduli e applicazioni''}
\maketitle
\section*{Introduction}

Recently there have been several results on the ample (or K\"{a}hler) cone of symplectic manifolds, generalising the well known case of $K3$ surfaces. This was mainly done through the study of stability conditions on derived categories of surfaces, which has been completed in several cases. This yielded a wall and chamber decomposition of the positive cone of moduli spaces of stable objects on those derived categories. Basically, there is a set of divisors whose orthogonals divide the positive cone in chambers and one such chamber is the ample cone. Moreover all other chambers can be reached by elementary birational transformations and reflections on contractible divisors. Such a result was obtained by Yoshioka \cite{yoshi} for smooth moduli spaces on abelian surfaces, by Bayer and Macr\`{i} \cite{bm} for smooth moduli spaces on $K3$ surfaces and by Meachan and Zhang \cite{mz} for some well behaved singular moduli spaces on $K3$ surfaces. It was implicitly expected that such a wall and chamber decomposition exists for all \issp manifolds, as several results pointed out, see for example \cite[Huybrechts]{huy_ann} for the general behaviour of birational maps between \issp manifolds and \cite[Hassett and Tschinkel]{ht2} for a detailed analysis of one of the known fourfolds. The wall and chamber decomposition obtained from the categorical point of view was also recently extended by Bayer, Hassett and Tschinkel \cite{bht} for all projective manifolds of \kntiposp and independently by the author \cite{me_kahl} also for the other deformation classes. In this paper we will exploit an interesting byproduct of \cite{me_kahl} to compute the monodromy group of \issp manifolds. The key idea is that parallel transport Hodge isometries preserve the wall and chamber decomposition of the positive cone. Knowing such a decomposition, it is possible to compute the group of isometries preserving it. This gives an upper bound for the group of parallel transport Hodge isometries and therefore also for the monodromy group. We will apply this idea to two of the known examples: deformations of generalised Kummers and of O'Grady's ten dimensional manifold \cite{ogr}. In the first case, thanks to previous computations by Markman \cite{marmeh}, we have a lower bound for the monodromy group and the upper bound we obtain coincides with it, therefore we obtain the expected monodromy group. In the second case, as stated in \cite[Conjecture 10.7]{mark}, the monodromy group is expected to coincide with all orientation preserving isometries. However we obtain an orientation preserving isometry which does not preserve a wall of the decomposition of the positive cone, therefore disproving the conjecture.    
\section*{Notations and preliminaries}
Let $X$ be an \issp manifold, \ie a simply connected \kahl manifold such that $H^{2,0}$ is generated by a symplectic form. For a detailed description of the properties of \issp manifolds we refer to Huybrechts survey \cite{huy}. We denote by $A_X$ the discriminant group of its second cohomology, \ie $H^2(X,\mathbb{Z})^\vee/H^2(X,\mathbb{Z})$, where the lattice structure is given by the Beauville-Bogomolov form $b(\,,\,)$. We also identify $H_2(X,\mathbb{Z})$ with a subset of $H^2(X,\mathbb{Q})$ by the natural embedding of the dual of a lattice $L$ inside $L\otimes\mathbb{Q}$. For any element $D$ of $H^2(X,\mathbb{Z})$, we denote with $div(D)$ a positive generator of the ideal $b(D,H^2(X,\mathbb{Z}))$. We remark that $D/div(D)$ is primitive in $H_2(X,\mathbb{Z})$.\\ We denote by $Mon^2(X)$, and call it the monodromy group, the group of parallel transport operators on the second cohomology of $X$, \ie the group of isometries of $H^2(X,\mathbb{Z})$ obtained by parallel transport and $Mon^2(X,Y)$ denotes parallel transport operators between two \issp manifolds. We remind the reader that $Mon^2(X)$ is a deformation invariant. $Hdg(X,Y)$ and $Hdg(X)$ denote respectively Hodge isometries between $H^2(X,\mathbb{Z})$ and $H^2(Y,\mathbb{Z})$ and Hodge isometries of $H^2(X,\mathbb{Z})$. The group $O^+(H^2(X,\mathbb{Z}))$ denotes orientation preserving isometries. For known results on monodromy group and related statements, we refer to Markman's survey \cite{mark}.\\
The positive cone $\mathcal{C}_X$ of an \issp manifold $X$ is the connected component of the cone of positive (with respect to the Beauville Bogomolov form) real $(1,1)$ classes containing the \kahl cone $\mathcal{K}_X$. The birational \kahl cone $\mathcal{BK}_X$ is the pullback of the \kahl cones of all \issp manifolds $X'$ birational to $X$. It is a disjoint union of convex cones and its closure coincides with the movable cone if $X$ is projective.\\
Given a symplectic surface $S$ and a polarisation $H$, $M_{v}(S,H)$ denotes the moduli space of $H$ stable sheaves on $S$ with Mukai vector $v\in H^{2*}(S)$. Here the Mukai vector of a sheaf $\mathcal{F}$ is $(rk(\mathcal{F}),c_1(\mathcal{F}),c_1(\mathcal{F})^2-c_2(\mathcal{F})+rk(\mathcal{F}))$ if $S$ is a $K3$ surface and $(rk(\mathcal{F}),c_1(\mathcal{F}),c_1(\mathcal{F})^2-c_2(\mathcal{F}))$ if $S$ is abelian. If $S$ is a $K3$ surface and $M_v(S,H)$ is smooth, then it is an \issp manifold and its second cohomology is isometric to $v^\perp\subset H^{2*}(S,\mathbb{Z})$. The same holds also for abelian surfaces after replacing $M_v(S,H)$ with its Albanese fibre $K_v(S,H)$. A similar construction can be extended to objects in the derived category of such surfaces and one obtains more deformations of the above mentioned manifolds, see \cite{bm} and \cite{yoshi}.\\
We call manifold of Kummer $n$ type any deformation of the Albanese fibre of $Hilb^{n+1}(S)$, where $S$ is an abelian surface. Such a fibre is $K_{(1,0,-n-1)}(S,H)$. These are \issp manifolds of dimension $2n$. The second example we consider is a symplectic resolution (and its smooth deformations) of a moduli space of sheaves on a $K3$ surface with Mukai vector $v=2w=(2,0,-2)$. Such example was constructed by O'Grady \cite{ogr}. Finally, we denote with $\Lambda_{24}$ the lattice isometric to $H^{2*}(K3,\mathbb{Z})$ and with $\Lambda_8$ the lattice $H^{2*}(A,\mathbb{Z})$ for an abelian surface $A$. They are respectively isometric to $U^4\oplus E_8(-1)^2$ and $U^4$, where $U$ is the hyperbolic lattice and $E_8$ is the unique positive even unimodular lattice of rank 8.

\section*{Acknowledgements}
I would like to thank E. Markman for his comments and suggestions, M. Wandel for reading a preliminary version of this paper and K.G. O'Grady and A. Rapagnetta for useful discussions. I am also grateful to the University of Lille, for their kind hospitality when the most relevant part of this work was carried out.  


\section{Decomposition of the positive cone and Monodromy}
The decomposition of the positive cone we are concerned with is a set of open chambers whose boundary walls are given by the orthogonal to the following divisors.
\begin{defn}
Let $X$ be an \issp manifold and let $D$ be a divisor on $X$. Then $D$ is called a wall divisor if $D^2<0$ and $h(D^\perp)\cap\mathcal{BK}_X=\emptyset$, for all parallel transport Hodge isometries $h$.
\end{defn}

The most important fact about such a decomposition is that one of the open chambers is the \kahl cone. Indeed, let $R$ be an extremal ray of the Mori cone of $X$ and let $D$ be a primitive divisor such that $D/div(D)=R$. Then in \cite[Lemma 1.4]{me_kahl} it is proven that $D$ is a wall divisor. The second fundamental property of wall divisors is that they are preserved under parallel transport Hodge isometries. 

\begin{thm}\cite[Theorem 1.3]{me_kahl}
Let $X$ and $Y$ be \issp manifolds and let $\mathcal{D}_X$ and $\mathcal{D}_Y$ be the sets of wall divisors of $X$ and $Y$. Then $Mon^2(X,Y) \cap Hdg^2(X,Y)$ sends $\mathcal{D}_X$ into $\mathcal{D}_Y$.
\end{thm}


The positive cone $\mathcal{C}_X$ of an \issp manifold is therefore divided in wall and chambers by the orthogonal to any element of $\mathcal{D}_X$. Using stability conditions on the derived category of a symplectic surface, this wall and chamber decomposition was determined for several \issp manifolds.

Let $M_v(S,H)$ be the moduli space of stable sheaves with primitive Mukai vector $v$ on the $K3$ surface $S$ with respect to a given $v$ generic polarization $H$. Let $D$ be a divisor of $M_v(S,H)$ with $D^2<0$. Let $T$ be the rank two hyperbolic primitive sublattice of $\Lambda_{24}$ containing $v$ and $D$. 
\begin{thm}\cite[Theorem 5.7 and 12.1]{bm}\label{thm:eman_muri}
Let $D,v$ and $T$ be as above, then $D$ is a wall divisor if and only if one of the following holds
\begin{itemize}
\item There exists $w\in T$ such that $w^2=-2$ and $0\leq (w,v)\leq v^2/2$.
\item There exists $w\in T$ such that $w^2\geq 0$ and $w^2<(w,v)\leq v^2/2$.
\end{itemize}
\end{thm} 

Similarly, a decomposition of the ample cone for the Albanese fibre of moduli spaces of stable objects on the derived category of an abelian surface has been obtained by Yoshioka.\\ Let $K_v(S,H)$ be the Albanese fibre of the moduli space of stable sheaves with primitive Mukai vector $v$ on the abelian surface $S$ with respect to a given $v$ generic polarization $H$. Let $D$ be a divisor of $K_v(S,H)$ with $D^2<0$. Let $T$ be the rank two hyperbolic primitive sublattice of $\Lambda_{8}$ containing $v$ and $D$.

\begin{thm}\cite[Proposition 1.2 and section 3.2]{yoshi}\label{thm:yoshi_muri}
Let $v,D$ and $T$ be as above, then $D$ is a wall divisor of $K_v(S,H)$ if and only if there exists $w\in T$ such that $w^2\geq 0$ and $w^2<(w,v)\leq v^2/2$.
\end{thm}  

Meachan and Zhang extended these results to special singular moduli spaces (and their symplectic resolutions), obtaining certain wall divisors for O'Grady's ten dimensional example.
Let $M_{2w}(S,H)$ be a $10$ dimensional singular moduli space of stable sheaves on a $K3$ surface $S$ with respect to a $2w$ generic polarisation $H$ and let $X$ be its symplectic resolution of singularities. Let $D$ be a divisor on $X$ obtained by pullback from a divisor $D'$ on $M_{2w}(S,H)$. Let $T$ be the rank two hyperbolic primitive sublattice of $\Lambda_{24}$ containing $w$ and $D'$.
\begin{thm}\cite[Theorem 5.3 and Theorem 5.4]{mz}\label{thm:meac_zhan_muri}
Keep notation as above. Then $D$ is a wall divisor on $X$ if and only if one of the following holds:
\begin{itemize}
\item There exists $s\in T$ such that $s^2=-2$ and $(s,w)=0$,
\item There exists $s\in T$ such that $s^2=-2$ and $(s,w)=1$.
\end{itemize} 
Let $D$ be a primitive generator of $\langle w,s\rangle\cap w^\perp$. Then $D$ is a wall divisor on $X$.
\end{thm}

We remark that these are the wall divisors for the moduli space $M_w(S,H)$, which sits naturally inside the singular locus of $M_{2w}(S,H)$.

\section{Kummer $n$ type manifolds}
In this section we compute the monodromy group of manifolds of Kummer $n$ type using two results. The first is Markman's computation of the intersection between the monodromy group of a Kummer $n$ type manifold $X$ with the group of isometries acting as $\pm1$ on $A_X$. The second is Yoshioka's result on the \kahl cone of generalised Kummers. We keep the same notation of \cite{marmeh}: let $\mathcal{W}(X)$ be the subgroup of $O^+(H^2(X,\mathbb{Z}))$ acting as $\pm 1$ on $A_X$. It is an order $2^a$ subroup of the group of orientation preserving isometries, where $a+1$ is the number of prime factors of $n+1$. Let $\chi$ denote the character corresponding to the action on $A_X$. Let $\mathcal{N}_X$ be the kernel of $det\circ \chi\,:\,\mathcal{W}_X\rightarrow \{\pm1\}$. 

\begin{prop}\cite[Corollary 4.8]{marmeh}
Let $X$ be a manifold of Kummer $n$ type. Then $Mon^2(X)\cap \mathcal{W}_X=\mathcal{N}_X$.
\end{prop}

\begin{prop}\cite{yoshi}
Let $X=K_{v}(S,H)$ for an abelian surface $S$. Then the movable cone of $X$ is cut out by all primitive isotropic elements $w\in H^{2*}(S,\mathbb{Z})$ of type $(1,1)$ such that $(v,w)=1$ or $2$.
\begin{proof}
Let $w$ be in the above set of isotropic classes. Then $v$ and $w$ generate a lattice $T$ which satisfies the hypothesis of \Ref{thm}{yoshi_muri}, so any generator of $v^\perp\cap T$ is a wall divisor. Notice that such a divisor $D$ is $v-(2n+2)w$ if $(v,w)=1$ and $v-(n+1)w$ in the second case. By \cite[Theorem 3.29]{yoshi}, there is a divisorial contraction associated to such elements and they cut out the movable cone.
\end{proof}  
\end{prop}
We remark that the two cases above belong to two different orbits of the isometry group of $H^2(X,\mathbb{Z})$, since the associated wall divisors have different divisibility. We will call them divisorial contractions of type I when $(v,w)=1$ and type II otherwise. In particular this type subdivision is preserved for all wall divisors in the same isometry orbit.

\begin{thm}\label{thm:mono_kum}
Let $X$ be a manifold of Kummer n type. Then $Mon^2(X)=\mathcal{N}_X$
\begin{proof}
To prove this statement, it is enough to prove that a wall preserving isometry must act as $\pm 1$ on the discriminant group $A_X$. Without loss of generality, we can do this computation on a generalised Kummer $X$. We denote by $2\delta$ the class of the exceptional divisor of the Hilbert-Chow morphism and $v=(1,0,-n-1)$. Notice that $\delta^\perp\subset H^2(X,\mathbb{Z})$ is unimodular. We let $s=\frac{v-\delta}{2n+2}\in H^{2*}(S,\mathbb{Z})$. The lattice $T=\langle v,s \rangle$ is isometric to $\left( \begin{array}{cc} 2n+2 & 1\\ 1 & 0 \end{array} \right)$, which means that $\delta$ is a divisorial contraction of type I. Let $g$ be an isometry of $H^2(X)$ sending $\delta$ to another wall divisor. We have $g(\delta)=k\delta+(2n+2)l$, where $l\in\delta^\perp$. As explained in \cite[Proposition 5.3]{yoshi}, either $t=\frac{v+g(\delta)}{2n+2}$ or $t'=\frac{v-g(\delta)}{2n+2}$ is an integer class of $H^{2*}(S)$. In both cases, we have that $t-l-s=\frac{(k-1)\delta}{2n+2}$ or $t-l+s=\frac{(k+1)\delta}{2n+2}$ are integer classes. This implies $k=\pm 1$ modulo $2n+2$, which means that the action of $g$ on $A_X$ is $\pm 1$. 
\end{proof}
\end{thm}

\begin{oss}
The above computation implies in particular that a manifold $X$ whose second cohomology is Hodge isometric to that of $K_v(S,H)$, for some $S$, $v$ and $H$, is itself a moduli space of stable objects on the derived category of the same abelian surface $S$. This was conjectured by Wandel and the author in \cite{mw}. The missing ingredient to extend the proof of \cite[Proposition 2.4]{mw} is precisely the above computation of the monodromy group.
\end{oss}

\section{O'Grady's ten dimensional manifold}
In this section we focus on the deformation class of the ten dimensional sporadic \issp manifold constructed by O'Grady \cite{ogr}. We cook up a peculiar example and we combine it with Meachan and Zhang's results to obtain a restriction on the monodromy group. In the following $S$ is a $K3$ surface, $v=2w$ a non primitive Mukai vector of square $8$, $H$ a $v$ generic polarisation on $S$ and $X$ is the symplectic resolution of $M_v(S,H)$. The manifold $X$ is deformation equivalent to O'Grady's ten dimensional manifold.
We will also use a natural embedding given by Perego and Rapagnetta \cite{pr} for the second cohomology of such manifolds.
\begin{thm}\cite[Theorem 1.7]{pr}
Keep notation as above. There exists a pure weight two Hodge structure on $H^2(M_v(S,H),\mathbb{Z})$, a Hodge isometry $v^\perp=w^\perp\rightarrow H^2(M_v(S,H),\mathbb{Z})$ and the pullback map $H^2(M_v(S,H),\mathbb{Z})\rightarrow H^2(X,\mathbb{Z})$ is injective.
\end{thm}

\begin{ex}
Let $S$ be a projective $K3$ surface with a symplectic automorphism $\varphi$ of order $3$ and let $X$ be a symplectic resolution of the moduli space of stable sheaves $M_v(S,H)$ for a $\varphi$ invariant ample class $H$. Here $v=2w$, $v^2=8$ and $\varphi(v)=v$. As remarked in \cite[Proposition 4.3]{mw}, $\varphi$ induces a symplectic automorphism $\widehat{\varphi}$ of $X$. Let $T_{\widehat{\varphi}}(X)$ be the invariant sublattice for the induced action of $\widehat{\varphi}$ on $H^2(X,\mathbb{Z})$ and let $S_{\widehat{\varphi}}(X)$ be its orthogonal. Again in \cite[Proposition 4.3]{mw} it is proven that $S_{\widehat{\varphi}}(X)\cong S_\varphi(S)$. Notice in particular that $S_{\widehat{\varphi}}(X)\subset H^{1,1}(X)$, since the symplectic form is invariant, and that there exists a $\widehat{\varphi}$ invariant \kahl class $\omega$. The lattice $S_{\varphi}(S)$ has been explicitely computed by Garbagnati and Sarti \cite{gs}. In particular, it contains an element $F$ of square $-10$, obviously orthogonal to the \kahl class $\omega$. Therefore $F$ is not a wall divisor on $X$.
\end{ex}

\begin{thm}
Let $X$ be a manifold deformation equivalent to O'Grady's ten dimensional \issp manifold. Then $Mon^2(X)$ is strictly smaller than $O^+(H^2(X,\mathbb{Z}))$.

\begin{proof}
Let $S$ be a very general $K3$ surface having a symplectic automorphism of order $3$ and an invariant polarisation $H$ of degree $2$.
Let $X$ be the resolution of singularities of a moduli space of stable objects in the derived category of $S$ with Mukai vector $v=2w=(2,0,-2)$. Let $s=(2,H,1)\in H^*(S,\mathbb{Z})$. We have $s^2=-2$ and $(s,w)=1.$ Let $T=\langle w,s\rangle$ and let $D$ be a generator of $T\cap v^\perp$. Then $D$ is a wall divisor on $X$, as follows from \Ref{thm}{meac_zhan_muri}. Notice that $D$ has divisibility $2$ when considered as an element of $H^2(M_v(S,H),\mathbb{Z})$ but it has divisibility $1$ inside $H^2(X,\mathbb{Z})$. Let now $F$ be (as in the previous example) a class of square $-10$ inside the coinvariant lattice of the induced order $3$ symplectic automorphism on $X$. By Eichler's criterion \cite[Lemma 3.5]{ghs}, there exists an orientation preserving isometry of $H^2(X,\mathbb{Z})$ sending $D$ to $F$. This implies that the group of isometries preserving wall divisors is strictly smaller than the group of isometries. Hence, also the monodromy group is strictly smaller than $O^+(H^2(X,\mathbb{Z}))$.  
\end{proof}

\end{thm}

As said before, this differs from what was previously expected, see \cite[Conjecture 10.7]{mark}.

\end{document}